\documentclass{amsart}
\usepackage[pdfauthor={Richard J. Mathar},pdfkeywords={Zernike Polynomial, Adaptive Optics},pdfsubject={Mathematics, Analysis},colorlinks=true,citecolor=blue]{hyperref}

\usepackage{amsmath}

\newtheorem{rem}{Remark}
\def\sgn{\mathop{\mathrm{sgn}}\nolimits}

\begin{document}

\title{A C++ Incarnation of Zernike Circle Functions}

\author{Richard J. Mathar} 
\email{mathar@mpia.de}
\urladdr{http://www.mpia.de/~mathar}
\address{R. J. Mathar, Hoeschstr. 7, 52372 Kreuzau, Germany}


\date{\today}

\begin{abstract}
An explicit C++ library is provided which deals with Zernike Functions over
the unit circle as the main subject. The implementation includes basic means
to evaluate the functions at points inside the unit circle and to convert
the radial and azimuthal parameters to Noll's index and vice versa.
Advanced methods allow to expand products of Zernike Functions
into sums of Zernike Functions, and to convert Zernike Functions to
polynomials over the two Cartesian coordinates and vice versa.
\end{abstract}

\maketitle 

\section{Basis Functions} 

The Zernike circle functions $Z$ are products of a radial polynomial $R$
and an azimuthal
sine or cosine function $A$:
\begin{equation}
Z_{n}^{(m)}(r,\varphi) = R_n^{m}(r)A_m(\varphi),
\end{equation}
where
\begin{equation}
r=\sqrt{x^2+y^2}, \quad 0\le r\le 1
\end{equation}
 is the distance to the origin of coordinates
and
\begin{equation}
x=r\cos\varphi;\quad y=r\sin\varphi
\end{equation}
 define the azimuth angle $\varphi$.
$n$ is the degree of $R$, a non-negative integer.
$m$ is one of $-n,-n+2,\ldots ,n-2,n$, such that $n-m$ is an even integer number.
The normalization chosen here is \cite{NollJOSA66}
\begin{equation}
\int_0^1 r R_n^{m}(r)R_{n'}^{m}(r) dr=\delta_{n,n'}; 
\quad \int_0^{2\pi} A_m(\varphi)A_{m'}(\varphi)d\varphi=\delta_{m,m'},
\end{equation}
\begin{equation}
\int_0^1 r dr \int_0^{2\pi} d\varphi Z_n^{m}(r,\varphi)Z_{n'}^{m'}(r,\varphi) =\delta_{n,n'}\delta_{m,m'},
\end{equation}
such that
\begin{equation}
R_{n}^{m}(r)\equiv \sqrt{2n+2}(-1)^{(n-|m|)/2} \binom{\frac{n+|m|}{2}}{\frac{n-|m|}{2}} r^{|m|}
{}_2F_1(-\frac{n-|m|}{2}, 1+\frac{n+|m|}{2}; 1+|m|;r^2),
\label{eq.R}
\end{equation}
and
\begin{equation}
A_m(\varphi) =\left\{\begin{array}{ll}
\cos(m\varphi)/\sqrt{\epsilon_m \pi}, & m\ge 0 ; \\
\sin(|m|\varphi)/\sqrt{\pi}, & m <0. \\
\end{array} \right. 
\end{equation}
where
\begin{equation}
\epsilon_m\equiv \left\{
\begin{array}{ll} 2, & m=0; \\ 1,& m\neq 0.\end{array}
\right.
\end{equation}

\begin{rem}
This is sligthly different from the (more common) notation in my earlier 
representation \cite{MatharArxiv0809,MatharSAJ179}; therefore some of the
equations are reproduced here where factors need to be replaced.
\end{rem}

The product expansions of the azimuth functions are \cite[\S II.E]{MatharArxiv0809}:
\begin{equation}
A_mA_{m'}=
\frac{1}{2\sqrt{\epsilon_m\epsilon_{m'}\pi}}\times
\left\{
\begin{array}{ll}
\sqrt{\epsilon_{|m-m'|}}A_{|m-m'|}+\sqrt{\epsilon_{m+m'}}A_{m+m'}
, & m\ge 0,m'\ge0;\\
A_{-|m+|m'||}
-\sgn(m-|m'|)A_{-|m-|m'||}
, & m\ge 0, m'<0;\\
\sqrt{\epsilon_{m-m'}}A_{|m-m'|}
- A_{|m+m'|}
, & m<0,m'<0.\\
\end{array}
\right.
\label{eq.Aprod}
\end{equation}
The product expansions of the radial functions are:
\begin{equation}
R_{n}^{m}(r)R_{n'}^{m'}(r) =\sum_{n''=m''}^{n+n'} g_{n,m,n',m',n'',m''}R_{n''}^{m''}(r),
\end{equation}
with projections \cite[\S II.E]{MatharArxiv0809}
\begin{multline}
g_{n_1,m_1,n_2,m_2,n_3,m_3}
=\int_0^1 r \prod_{j=1}^3 R_{n_j}^{m_j}(r) dr
\\
=\sqrt{8\prod_{j=1}^3(n_j+1)}
\sum_{s_1=0}^{(n_1-|m_1|)/2}
\sum_{s_2=0}^{(n_2-|m_2|)/2}
\sum_{s_3=0}^{(n_3-|m_3|)/2}
\frac{1}{2+n_1+n_2+n_3-2(s_1+s_2+s_3)}
\\
\times
\prod_{j=1}^3
(-)^{s_j}\binom{n_j-s_j}{s_j}\binom{n_j-2s_j}{(n_j-|m_j|)/2-s_j}
.
\label{eq.g}
\end{multline}

The class functions which are described in the next section
answer the following questions:
\begin{enumerate}
\item
Given the numerical coefficients $c_{n,m}$ and $c_{n',m'}$,
what are the values of the linearization coefficients $c_{n'',m''}$
in
$\sum_{n,m} c_{n,m}Z_n^{(m)} \sum_{n',m'} c_{n'm'}Z_{n'}^{(m')}$ =$\sum_{n'',m''} c_{n'',m''}Z_{n''}^{(m'')}$?
\item
Given the numerical coefficients $c_{p,q}$,
what are the values of the coefficients $c_{n,m}$
in
$\sum_{p,q} c_{p,q}x^py^q=\sum_{n,m} c_{n,m}Z_{n}^{(m)}$?
\item
Vice versa, given the numerical coefficients $c_{n,m}$,
what are the values of the coefficients $c_{p,q}$
in
$ \sum_{n,m} c_{n,m}Z_{n}^{(m)}
=
\sum_{p,q} c_{p,q}x^py^q$?
\end{enumerate}

\section{Implementation}
\subsection{Points, Locations}
Points in the unit circle and in the unit sphere
are represented by their 2 to 3 Cartesian coordinates and implemented
in the classes \texttt{Point2D} and \texttt{Point3D}.
The constructors expect the 2 or 3 Cartesian coordinates that fix
the point in space. The \texttt{Point2D} object has a trivial method
to extract the radial distance $r$ and the azimuth angle
$\varphi$ of the circular coordinates from the $x$ and $y$ components.

\subsection{Polynomials}
A single term $cx^j$ of a univariate polynomial is represented
by an object of the class \texttt{Monomial1D}.
A single term $cx^py^q$ of a bivariate polynomial is represented
by an object of the class \texttt{Monomial2D}.
A single term $cx^py^qz^r$ of a trivariate polynomial is represented
by an object of the class \texttt{Monomial3D}.
The constructors accept the exponents $j$, the exponents $p$ and $q$
and the exponents $p$, $q$ and $r$ respectively. The coefficient $c$
is optional and set to unity if missing in the constructor.
The useful operations within these classes are multiplication with
or division through a constant value (which scales the coefficient),
and multiplication with a single term of the same type (which essentially
adds the exponents of the two factors). These operations are implemented
by overloaded multiplication and division operators.
There is a common method \texttt{at()} which evaluates the term
given a value $x$ on the line or a position specified by a \texttt{Point2D}
or \texttt{Point3D} in the plane or in three-dimensional space.

Univariate polynomials $\sum_j c_j x^j$, bivariate polynomials $\sum c_{p,q}x^py^q$
and trivariate polynomials $\sum c_{p,q,r} x^py^qz^r$ are represented
as vectors of the monomial objects in the classes \texttt{Polynomial1D},
\texttt{Polynomial2D} and \texttt{Polynomial3D}. Addition, subtraction
and multiplications within each of them are closed operations and
implemented by overloaded addition, subtraction and multiplication operators.
Adding new terms is supported with overloaded \texttt{+=} operators
that accept a single monomial object or another polynomial object
of the same dimension.

There is a common method \texttt{at()} which evaluates these polynomials
given a value $x$ on the line or a position specified by a \texttt{Point2D}
or \texttt{Point3D}. The function sums up the single-term components.

A special variant of the univariate polynomials are the terminating
Gaussian Hypergeometric Functions $_2F_1(a,b;c;z)$ with integer parameter
$a\le 0$, two auxiliary parameters $b$ and $c$, as a function of the
variable $z$. These are implemented as a class \texttt{Hypergeom21}
derived from \texttt{Polynomial1D}.

\subsection{Zernike Radial Function}
A radial function $R_n^m(r)$ is special case of the Hypergeometric
Function constructed as in Equation (\ref{eq.R}) given the parameter
$n\ge 0$ and the parameter $m$ (the latter being referred to only as $|m|$).
It is represented by an object of the class \texttt{ZernikeRadi},
a subclass of \texttt{Hypergeom21}.
A method in the class computes the $g$-factors of Equation (\ref{eq.g})
to support expansion of products of $R$-functions in other $R$-functions.

The implementation is wider than actually needed here:
an additional argument (which defaults to 2) supports use of
the radial function in $D\ge 2$ dimensions, where $D$ appears
on the right hand sides of Eqs.\ (\ref{eq.R}) and (\ref{eq.g}) \cite{MatharSAJ179}.

\subsection{Zernike Azimuthal Function}
A sine or cosine term of the form $cA_m$ is represented by
an object of the class \texttt{ZernikeCircAzi}, which is constructed
given the signed integer parameter $m$ and an optional prefactor $c$---which
is set to unity of missing. Scaling of the prefactor by a constant
is implemented by overloaded multiplication and division.
Evaluation at some explicit angle $\varphi$ happens by calling
the \texttt{at()} member function with an argument $\varphi$ in units of radians.

A collection (arithmetic sum) of the form $\sum_m c_m A_m(\varphi)$ is represented
by an object of the class \texttt{ZernikeCircAziVec}, where each
component is stored as an element of the \texttt{ZernikeCircAzi} type.
Adding new terms is achieved by using the overloaded \texttt{+=}
operator.
Scaling all terms (prefactors) at the same time is supported
by overloaded \texttt{*=} operation with a constant.
An overloaded multiplication implements the arithmetic
multiplication by means of Equation (\ref{eq.Aprod}).

\subsection{Zernike Circle Function}
A Zernike Circle Function $cZ_n^{(m)}$ is represented by an object
of the class \texttt{ZernikeCirc} which is initialized by the
radial parameter $n$, the azimuthal parameter $m$ and an optional
prefactor $c$---which is set to unity if missing.
Violation of the two constraints on $m$ (parity and range)
are silently caught by setting the prefactor to $c=0$.
Instead of the two parameters $n$ and $m$, Noll's index $j\ge 1$
can also be used to define an object of \texttt{ZernikeCirc}.

Scaling of the prefactor $c$ is supported by overloaded
\texttt{*=} and \texttt{/=} operators.

Evaluation of the function at some point in the unit circle
is done by calling the \texttt{at()} member function with an
argument/location specified by a \texttt{Point2D} object. The implementation
constructs the factors $R_n^m$ and $A_m$ with objects of the \texttt{ZernikeRadi}
and \texttt{ZernikeCircAzi} types and multiplies their values.
A design decision is to set the values to zero outside the unit circle, $r> 1$(!).

An arithmetic sum of the form $\sum c_{n,m} Z_n^{(m)}$ is represented
by an object of the \texttt{ZernikeCircVec} class, which represents
the components as members of a vector of \texttt{ZernikeCirc}'s.
Adding terms to the sum is supported by overloaded \texttt{+=} operators
for additional \texttt{ZernikeCirc} or \texttt{ZernikeCircVec} terms.
Scaling all the coefficients $c_{n,m}$ is achieved by overloaded
\texttt{*=} and \texttt{/=} operations with constant arguments.

The most valuable functions of the implementation are:
\begin{enumerate}
\item
The arithmetic product of two Zernike expansions is implemented
by the overloaded operator \texttt{*=} which works with two
factors of the \texttt{ZernikeCirc} or \texttt{ZernikeCircVec} type
and produces a \texttt{ZernikeCircVec} object. This is an incarnation
of \cite[\S II.E]{MatharArxiv0809}; it starts with the product
expansion of the $A$-terms provided the \texttt{ZernikeCircAziVec}
class and collects the product expansion of the $R$-terms by
calling the \texttt{g} member functions in the \texttt{ZernikeRadi}
class.
\item
The conversion of a polynomial $\sum cx^py^q$ form into
a Zernike expansion is supported by constructors in the \texttt{ZernikeCircVec}
class that accept arguments of the \texttt{Monomial2D} or
\texttt{Polynomial2D} type. This implements \cite[\S II.C]{MatharArxiv0809}.
\item
The conversion of a Zernike expansion $\sum c_{n,m}Z_n^{(m)}$
into a polynomial $\sum cx^py^q$
is supported by constructors in the \texttt{Polynomial2D}
class that accept arguments of the \texttt{ZernikeCirc} or
\texttt{ZernikeCircVec} type. This implements \cite[\S II.D]{MatharArxiv0809}.
\end{enumerate}

\subsection{Polynomial Fit}
If the GNU Scientific Library is available \cite{GSL}, an additional class \texttt{PowFit2D}
is introduced, which is fed with (constructed from) a list of Cartesian $x$ and $y$ coordinates and
function values $f(x,y)$ at these points. There are two simple formats of entering
the data into the constructor, one that reads the $x,y,f$ triples from an ASCII file,
the other providing them as a list of \texttt{Point3D} points (interpreting the third coordinate as $f$).

The major member function is the \texttt{fit()} function which constructs
the ordinary least squares fit through the $f$ values up to some total order of the
fitting polynomial $\sum \alpha x^py^q$. The limiting order $p+q$ is an argument of the \texttt{fit()}.
This is implemented as another constructor for 
the \texttt{Polynomial2D} class that admits a \texttt{PowFit2D} object
as its argument. Fitting of scattered data over the unit circle
to Zernike Circle Polynomials for some upper limit of the index $n$
is then a matter of converting the fitting \texttt{Polynomial2D}
to a \texttt{ZernikeCircVec} object in a final step.

Note that this is a  demonstration of library usage, but the strategy is
inefficient. In practise, the Zernike coefficients would be obtained directly by solving the
linear algebra in the $(r,\varphi)$ coordinates.

\bibliographystyle{amsplain}
\bibliography{all}

\appendix 
 
\section{Installation}
\subsection{Compilation}
The source code of roughly 4300 lines
is available in the \texttt{anc} directory, licensed under the GNU General Public License.

It is compiled with the GNU autotools \cite{GNUauto} via
\begin{verbatim}
autoreconf -i
./configure --prefix=$HOME
make
make install
\end{verbatim}
This bundles the classes in a library \texttt{libZernikeCirc.a} and
compiles the test routine \texttt{tstZernikeCirc}. It also searches
for the \texttt{gsl} and \texttt{gslcblas} libraries and includes
\texttt{PowFit2D} if these are found.

If the autotools are not
available, compilation in the conventional style of
\begin{verbatim}
g++ -c -O2 [A-Z]*.cxx
ld -i -o libZernikeCirc.a *.o
g++ -o tstZernikeCirc tstZernikeCirc.cxx -L. -l ZernikeCirc
\end{verbatim}
is an alternative. Definition of the preprocessor variables
\texttt{HAVE\_GSL\_GSL\_SF\_H} and \texttt{HAVE\_GSL\_GSL\_LINALG\_H}
and adding the flags \texttt{-lgsl -lgslcblas} must be done manually then,
if applicable.

If \texttt{doxygen} is available, the API documentation can be constructed in the \texttt{html} directory
with
\begin{verbatim}
make doc
firefox html/index.html
\end{verbatim}

\subsection{Numerical Tests}
The test program can be run with
\begin{verbatim}
tstZernikeCirc
\end{verbatim}
The test suite contains
\begin{enumerate}
\item
a table of mappings of the two parameters $(n,m)$ onto Noll's index $j$
to test the \texttt{nollIdx} member function of \texttt{ZernikeCirc};
\item
a table of mappings of Noll's index $j$ to the two parameters $(n,m)$
to test basically the inverse functionality in the constructor of \texttt{ZernikeCirc};
\item
a double loop over pairs of objects of \texttt{ZernikeCirc}
to test that their products \texttt{ZernikeCircVec} have the same
value at some \texttt{Point2D} as expected from the product of the
individual values;
\item
a loop over various terms of the $cx^py^q$ format to test that
the conversion of a \texttt{Monomial2D} object
into a \texttt{ZernikeCircVec} object
keeps its
value for some points \texttt{Point2D} scattered in the unit circle;
\item
a loop over various terms of the $cZ_n^{(m)}$ format to test
that the conversion of a \texttt{ZernikeCirc} object
into a \texttt{Polynomial2D} object
keeps its
value for some  points \texttt{Point2D} scattered in the unit circle.
\end{enumerate}

\end{document}